\documentclass[titlepage,11pt]{article}
\usepackage{latexsym}
\usepackage{amsfonts}
\usepackage{amsmath}
\usepackage{amsthm,amssymb,bm} 
\newcommand\blackslug{\hbox{\hskip 1pt \vrule width 4pt height 8pt depth 1.5pt
        \hskip 1pt}}
\newcommand\bbox{\hfill \quad \blackslug \bigbreak}

\def\E{{\mathbb E}}
\def\P{{\mathbb P}}

\title{Excluding pairs of graphs}
\author{Maria Chudnovsky\thanks{Supported by NSF grants DMS-1001091 and IIS-1117631.}\\
Columbia University, New York, NY 10027, USA
\\
\\
Alex Scott\\
Mathematical Institute, University of Oxford, 24-29 St Giles', Oxford OX1 3LB, 
UK	
\\
\\
Paul Seymour\thanks{Supported by ONR grant N00014-10-1-0680 and NSF grant DMS-0901075.}\\
Princeton University, Princeton, NJ 08544, USA}

\date{September 3, 2012; revised \today}

\newtheorem{theorem}{}[section]

\newcommand{\Proof}{\noindent{\bf Proof.}\ \ }

\begin{document}
\maketitle
\begin{abstract}

For a graph $G$ and a set of graphs $\mathcal{H}$, we say that $G$ is 
{\em $\mathcal{H}$-free} if no induced 
subgraph of $G$ is isomorphic to a member of $\mathcal{H}$. Given an integer 
$P>0$, a graph $G$, 
and a set of graphs $\mathcal{F}$, we say that $G$ 
{\em admits an $(\mathcal{F},P)$-partition} if the vertex set of $G$ 
can be partitioned into $P$ subsets  $X_1, \ldots, X_P$, so that for every 
$i \in \{1, \ldots, P\}$, either $|X_i|=1$, or 
the subgraph of $G$ induced by $X_i$ is 
$\{F\}$-free  for some $F \in \mathcal{F}$.

Our first result is the following. For every pair $(H,J)$ of graphs such that  
$H$ is  the disjoint union of two  graphs $H_1$ and $H_2$, and the complement 
$J^c$ of $J$ is 
the disjoint union of two  graphs $J_1^c$ and $J_2^c$, there exists
an integer $P>0$ such that every $\{H,J\}$-free graph has an
$(\{H_1,H_2,J_1,J_2\},P)$-partition.
Using a similar idea we also give a short proof of one of the results of
\cite{heroes}.

A {\em cograph} is a graph obtained from single vertices by repeatedly taking
disjoint unions and disjoint unions in the complement. For every cograph there
is a parameter measuring its complexity, called its {\em height}.  
Given a graph $G$ and a pair of graphs $H_1,H_2$, we say that $G$
is {\em $\{H_1,H_2\}$-split} if $V(G)=X_1 \cup X_2$, where the subgraph of $G$ 
induced by $X_i$ is $\{H_i\}$-free for every $i \in \{1,2\}$.
Our second result is that for every integer $k>0$ and  pair $\{H,J\}$ of 
cographs each of height  $k+1$, where neither of $H,J^c$ is
connected, there exists a pair of cographs $(\tilde{H},\tilde{J})$,
each of height $k$, where neither of $\tilde{H}^c,\tilde{J}$
is connected, such that every $\{H,J\}$-free graph  is 
$\{\tilde{H},\tilde{J}\}$-split.

Our final result is  a construction showing that if $\{H,J\}$ are graphs
each with at least one edge, then for every pair  of integers $r,k$ there
exists a graph $G$ such that every $r$-vertex induced subgraph of $G$ is 
$\{H,J\}$-split, but $G$ does not admits an $(\{H,J\},k)$-partition.
\end{abstract}

\section{Introduction}

All graphs in this paper are finite and simple. Let $G$ be a graph.
For $X \subseteq V(G)$, we denote by $G|X$ the subgraph of $G$ induced
by $X$. The complement of $G$, denoted by $G^c$, is  the graph with vertex set 
$V(G)$ such that two vertices are adjacent in $G$ if and only if they are 
non-adjacent in $G^c$. A {\em clique} in $G$ is a set of vertices all
pairwise adjacent; and a {\em stable set} in $G$ is a set of vertices
all pairwise non-adjacent.

We denote by $K_n$ the complete graph on $n$
vertices, and by $S_n$ the complement of $K_n$.
For graphs $H$ and $G$ we say that $G$ {\em contains} $H$ if
some induced subgraph of $G$ is isomorphic to $H$.
Let $\mathcal{F}$ be a set of graphs. We say that $G$ is  
{\em $\mathcal{F}$-free} if $G$ contains no member of $\mathcal{F}$.
If $|\mathcal{F}|=1$, say $\mathcal{F}=\{F\}$, we write ``$G$ is $F$-free''
instead of ``$G$ is $\{F\}$-free''. For a pair of graphs
$\{H_1, H_2\}$, we say that $G$ is {\em $\{H_1,H_2\}$-split} if 
$V(G)=X_1 \cup X_2$, and $G|X_i$ is $H_i$-free for $i=1,2$. 
We remind the reader that a {\em split graph} is a graph whose vertex set
can be partitioned into a clique and a stable set; thus in our language
split graphs are precisely the graphs that are $\{K_2,S_2\}$-split.
Given an integer  $P>0$, we say that $G$
{\em admits an $(\mathcal{F},P)$-partition} if
$V(X)=X_1 \cup \ldots \cup X_P$ such that for every $i \in \{1,\ldots, P\}$,
either $|X_i|=1$ or $G|X_i$  is $\{F\}$-free  for some $F \in \mathcal{F}$.
Please note that the first alternative in the definition of an  
$(\mathcal{F},P)$-partition (the condition that $|X_i|=1$) is only
necessary when no graph in $\mathcal{F}$ has more than one vertex.

In \cite{undirectedkillers} two of us proved the following:

\begin{theorem}
\label{multipartite}
For every pair of graphs $(H,J)$ such that $H^c$ and $J$ are complete 
multipartite, there exist integers $k,P>0$ such that every $\{H,J\}$-free 
graph admits a $(\{K_k,S_k\}, P)$-partition.
\end{theorem}
\noindent and its immediate corollary:
\begin{theorem}
\label{Ksplit}
For every pair of graphs $(H,J)$ such that $H^c$ and $J$ are complete 
multipartite, there exists an integer $k>0$ such that every $\{H,J\}$-free 
graph is $\{K_k,S_k\}$-split.
\end{theorem}

The first goal of this paper is to generalize \ref{multipartite}. Let
$H$ be a graph. A {\em component} of $H$ is a maximal connected subgraph of 
$H$. A graph is {\em anticonnected} if its complement is connected. An 
{\em anticomponent} of $H$ is a maximal anticonnected induced subgraph of $H$.
We denote by $c(H)$ the set of components of $H$, and by $ac(H)$
the set of anticomponents of $H$. We remark that for every non-null graph $G$,
at least one of $c(G)$ or $ac(G)$ equals $\{G\}$. We prove the following
generalization of \ref{multipartite} (please note that \ref{disconnected} is 
trivial unless $H$ is not connected and $G$ is not anticonnected):

\begin{theorem}
\label{disconnected}
For every pair of graphs $(H,J)$ there exists an integer $P$ such that
every $\{H,J\}$-free  graph admits a $(c(H) \cup ac(J), P)$-partition.
\end{theorem}

Please note that applying \ref{disconnected} with $J$ a complete
graph and $H$ a graph with no edges gives  Ramsey's famous
theorem. Using ideas similar to those of our proof of \ref{disconnected}, we 
also give a short proof of one of the results of \cite{heroes}.

Next let us generalize the notion of a complete multipartite graph. A 
{\em cograph} is a graph obtained from $1$-vertex graphs by repeatedly taking 
disjoint unions and disjoint unions in the complement. In particular, $G$ is
either not connected or not anticonnected for every cograph $G$ with at 
least two vertices,
and therefore for every cograph $G$ with at least two vertices, exactly one of $G,G^c$ is connected. We
recursively define a parameter, called the {\em height}  of a cograph, that 
measures its complexity, as follows. The height of a one vertex cograph is zero. If $G$
is a cograph that is not connected, let $m$ be the maximum
height of a component of $G$; then the height of $G$ is $m+1$. 
If $G$ is a cograph that is not anticonnected, let $m$ be the maximum
height of an anticomponent of $G$; then the height of $G$ is $m+1$. 
We denote the height of $G$ by $h(G)$.


We use \ref{disconnected} to prove the following:

\begin{theorem}
\label{cographsplit}
Let $k > 0$ be an integer, and let  $H$ and $J$ be cographs, each of height
$k+1$, such that $H$ is anticonnected, and $J$ is connected. Then there
exist cographs $\tilde{H}$ and $\tilde{J}$, each of height $k$, such that
 $\tilde{H}$ is connected, and $\tilde{J}$ is anticonnected, and 
every $\{H,J\}$-free  graph is $\{\tilde{H},\tilde{J}\}$-split.
\end{theorem}

The proof of \ref{multipartite} in \cite{undirectedkillers} relies on
the following lemma:
\begin{theorem}
\label{bounded1}
Let $p>0$ be an  integer. There exist  an integer $r>0$ such that
for every graph $G$, if every induced subgraph of $G$  with at most $r$ 
vertices is  $\{K_p, S_p\}$-split, then $G$ is  
$\{K_p, S_p\}$-split.
\end{theorem}

Here is a weaker statement that would still imply the results of 
\cite{undirectedkillers}:

\begin{theorem}
\label{bounded}
Let $p>0$ be an  integer. There exist  integers $r,k>0$ such that
for every graph $G$, if every induced subgraph of $G$  with at most $r$ 
vertices is  $\{K_p, S_p\}$-split, then $G$ admits  a 
$(\{K_p, S_p\},k)$-partition.
\end{theorem}

Originally we hoped that \ref{cographsplit} could be proved along the same 
lines, and that a result similar to \ref{bounded} might exist when the pairs 
$\{K_p,S_p\}$ were replaced by  pairs of more general graphs.  However, this 
turns out not to be the case, because of the following:

\begin{theorem}
\label{construction}
Let $H$, $J$ be graphs each with at least one edge.  Then for any choice of integers $r,k$ there is a graph $G$ such that 
\begin{itemize}
\item for every $S \subseteq V(G)$ with $|S| \leq r$, the graph
$G|S$ is $\{H,J\}$-split, and
\item $G$ has no $(\{H,J\}, k)$-partition.
\end{itemize}
\end{theorem}

By taking complements, the conclusion of \ref{construction} also holds
if each of $H$ and $J$ has a non-edge. Thus \ref{bounded} is in
a sense the strongest result of this form possible.

This paper is organized as follows. In Section~2 we prove~\ref{disconnected}.
In Section~3 we use \ref{disconnected} to prove~\ref{cographsplit}. In 
Section~4 we reprove a result of \cite{heroes}.
Finally, Section~5 contains the proof of \ref{construction}.

\section{The proof of \ref{disconnected}}

The goal of this section is to prove~\ref{disconnected}. 
Let us start with some definitions. Let $G$ be a graph. For disjoint
$X,Y \subseteq V(G)$, we say that $X$ is {\em complete} ({\em anticomplete})
to $Y$ if every vertex of $X$ is adjacent (non-adjacent) to every vertex of
$Y$. If $|X|=1$, say $X=\{x\}$, we say ``$x$ is complete (anticomplete)
to $Y$'' instead of ``$\{x\}$ is complete (anticomplete)
to $Y$''.

For a graph $G$ and a set 
$X \subseteq V(G)$, we denote by $G \setminus X$ the graph 
$G|(V(G) \setminus X)$. If $|X|=1$, say $X=\{v\}$, we write $G \setminus v$
instead of $G \setminus \{v\}$. Let $S$, $T$ be  induced subgraphs of a graph
$H^*$. We say that  $S$ is an {\em $H^*$-extension} of $T$ if $T=S \setminus u$ 
for some  $u \in V(S)$.

Clearly, \ref{disconnected} follows from repeated applications of the following:

\begin{theorem}
\label{twographs}
Let $H$, $J$, $H_1$, $H_2$, $J_1$ and $J_2$ be non-null graphs such that
$H$ is the disjoint union of $H_1$ and $H_2$, and
$J^c$ is the disjoint union of $J_1^c$ and $J_2^c$.
Let $m=\max(|V(H_1)|,|V(H_2)|,|V(J_1)|,|V(J_2)|)$, and $\mathcal{F} = \{H_1,H_2,J_1,J_2\}$.
Then every $\{H,J\}$-free graph $G$ admits an 
$(\mathcal{F}, 2(m+1)^m)$-partition.
\end{theorem}

\Proof We may assume without loss of generality that $|V(H_1)|=m$.
We may also assume that $G$ is not $H_1$-free, for otherwise the theorem holds.
Choose some vertex $h^*\in V(H_2)$, and let $H^*$ be the subgraph of $H$ induced on $V(H_1)\cup \{h^*\}$.
For an induced subgraph $T$ of $H^*$, a  {\em $T$-piece} is an isomorphism
$g$ from $T$ to an induced subgraph of $G$, that we denote by $g(T)$.
For a $T$-piece $g$, and an
$H^*$-extension $S$ of $T$, let $V(S)\setminus V(T) = \{s\}$; we say that $v \in V(G) \setminus V(g(T))$ {\em $g$-corresponds to} $S$ if
the map sending $x$ to $g(x)$ for each $x\in V(T)$, and sending $s$ to $v$, is an $S$-piece.
We denote by $Y_S(g)$ the set of all vertices in $G$ that $g$-correspond to $S$.
Let $Y(g)$ be
the union of the sets $Y_S(g)$. Thus $Y(g)$ is the set of all 
vertices in $V(G)$ that $g$-correspond to an $H^*$-extension of 
$T$.
\\
\\
(1) {\em If $g$ is an $H_1$-piece, then  
$Y(g)$ is anticomplete to  $V(g(H_1))$, and the graph $G|Y(g)$ is $H_2$-free.}
\\
\\
Since $V(H_1)$ is anticomplete to $V(H_2)$ in $H$, and in particular anticomplete to $h^*$, it follows that
$Y(g)$ is anticomplete to $V(g(H_1))$. But then, since 
$G$ is $H$-free, it follows 
that $G|Y(g)$ is $H_2$-free. This proves~(1).
\\
\\
Let 
$$\phi(t)= 
\begin{cases} 2(m+1)^{m-t} &\mbox{for } 0 \leq t \leq m-1\\
              1 &\mbox{for } t = m.
\end{cases}$$
\\
\\
(2) {\em Let $T$ be an induced subgraph of $H_1$, and let $g$
be a $T$-piece. Write $t=|V(T)|$. Then for some $H^*$-extension $S$ of $T$,
the graph  $G|Y_S(g)$ admits an $(\mathcal{F}, \phi(t))$-partition.}
\\
\\
The proof is by induction on $m-t$. 
If $t=m$, then (as $|H_1|=m$) it follows from~(1) that $G|Y(g)$ admits an $(\mathcal{F}, 1)$-partition,
so we may assume that $t \leq m-1$.

Choose an $H^*$-extension $S$ of $T$, with $S$ an induced subgraph of $H_1$.
If the graph  $G|Y_S(g)$ admits an $(\mathcal{F}, \phi(t))$-partition then we are done;
so we may assume that for every partition
of $Y_S(g)$ into at most $\phi(t)$ classes, some class contains each of
$H_1$, $H_2$, $J_1$, $J_2$.

Let $V(S)\setminus V(T) = \{s\}$.
Let $v \in Y_S(g)$, and let $h$ be the $S$-piece mapping $s$ to $v$ and mapping $x$ to $g(x)$ for each $x\in V(T)$.
Inductively, there exists an $H^*$-extension $Q$ of $S$
such that $G|Y_Q(h)$ admits an  
$(\mathcal{F}, \phi(t+1))$-partition.
Let $V(Q)\setminus V(S) = \{r\}$, and let $R = Q\setminus s$.
Thus $R$ is an $H^*$-extension of $T$ different from $S$. (Nevertheless, possibly $Y_R(g) = Y_S(g)$, if there is an isomorphism from $S$ to $R$
fixing $T$ pointwise.)
We say that $v$ is {\em of type $R$}.
From the definition of $Y_Q(h)$, it follows that
$Y_Q(h) \subseteq Y_R(g)$,
and either
\begin{itemize}
\item $r,s$ are non-adjacent in $H$, and $Y_Q(h)$ is the set of vertices in $Y_R(g)$ that
are different from and non-adjacent to $v$
or  
\item $r,s$ are adjacent in $H$, and $Y_Q(h)$ is the set of 
vertices in $Y_R(g)$ that are different from and adjacent to $v$.
\end{itemize}

For each $H^*$-extension $R$ of $T$ different from $S$, let $Z_{R}$ be the set of vertices in $Y_S(g)$ that are of type $R$.
Thus the sets $Z_{R}$ have union $Y_S(g)$.
Since $|V(H^*)|-|V(T)| = m+1-t$ and $R$ must be different from $S$, it follows that
are at most $m-t$ different types of vertices in $Y_S(g)$.
Since 
$$m-t \leq 2(m+1)^{m-t}=\phi(t),$$
there is an $H^*$-extension $R$ of $T$ different from $S$, such that
$G|Z_{R}$  contains each of
$H_1$, $H_2$, $J_1$, $J_2$. Write $A=Y_R(g)$.
Let $V(R)\setminus V(T) = \{r\}$.

Assume first that $r,s$ are non-adjacent,
and choose $B \subseteq Z_{R}$ such that
$G|B$ is isomorphic to $J_1$. For each $b \in B$, let $A_b$ be the set
of vertices of $A \setminus B$ that are non-adjacent to $b$.
Let $A_0$ be the set of vertices of $A \setminus B$ that are complete
to $B$. Then $A \subseteq B \cup A_0 \cup \bigcup_{b \in B}A_b$.

Let $b \in B$; since $b \in Z_{R}$, since
$A_b \subseteq Y_R(g)$, and since $A_b$ is anticomplete to $b$,
it follows from the definition of $Z_{R}$ that $G|A_b$ admits an 
$(\mathcal{F}, \phi(t+1))$-partition.
Since $G$ is $J$-free, it follows that $G|A_0$ is $J_2$-free.
Since $|B| \leq m$, this  implies that $G|A$ admits an 
$(\mathcal{F},K)$-partition, where 
$$K= m \phi(t+1)+ m +1 \leq 2(m+1)^{m-t}=\phi(t),$$ 
(even when $t=m-1$), and so (2) holds.

Now we assume that $r,s$ are adjacent.
Choose $B \subseteq Z_{R}$ such that
$G|B$ is isomorphic to $H_1$. For $b \in B$, let $A_b$ be the set
of vertices of $A \setminus B$ that are adjacent to $b$.
Let $A_0$ be the set of vertices of $A \setminus B$ that are anticomplete
to $B$. Then $A \subseteq B \cup A_0 \cup \bigcup_{b \in B}A_b$.

Let $b \in B$;  since $b \in Z_{R}$, since
$A_b \subseteq Y_R(g)$, and since $A_b$ is complete to $b$, 
it follows from the definition of $Z_{R}$ that $G|A_b$ admits an 
$(\mathcal{F}, \phi(t+1))$-partition.
Since $G$ is $H$-free, it follows that $G|A_0$ is $H_2$-free.
Since $|B| \leq m$, this  implies that $G|A$ admits an  
$(\mathcal{F},K)$-partition, where 
$$K = m \phi(t+1) + m +1 \leq 2(m+1)^{m-t}=\phi(t).$$ 
This proves~(2).
\\
\\
Now let $T$ be the null graph, and $g$ the isomorphism from $T$ into $G$. Then $T$ is an 
induced subgraph of $H_1$, and $g$ is a $T$-piece. Also, 
$Y_S(g) = V(G)$ for every $H^*$-extension $S$ of $T$.
But then, by~(2), $G$ admits an
$(\mathcal{F}, 2(m+1)^m)$-partition. This proves~\ref{twographs}.~\bbox

\section{Tournament heroes}
A {\em tournament} is a digraph such that for every two distinct vertices $u,v$
there is exactly one edge with ends $\{u,v\}$ (so, either the edge $uv$ or $vu$ but not both). Let $G$ be a tournament. If $uv$ is an edge of $G$ we say that
$u$ is {\em adjacent to} $v$, and $v$ is {\em adjacent from} $u$.
For $X \subseteq V(G)$, we denote by $G|X$ the subtournament
of $G$ induced by  $X$. We write $G \setminus X$ to mean 
$G|(V(G) \setminus X)$;  and if $|X|=1$, say $X=\{x\}$, 
we write $G\setminus x$ instead of $G \setminus \{x\}$.

If $X$ and $Y$ are two 
disjoint subsets of $V(G)$, we say that $X$ is {\em complete to} $Y$, 
and $Y$ is {\em complete from} $X$, if every vertex in $X$ is adjacent to every
vertex in $Y$; if $|X|=1$, say $X=\{x\}$, we say that $x$ is complete to $Y$, 
and $Y$ is complete from $x$. If $H$ is a  tournament, we say $G$ 
{\em contains} $H$ if $H$ is isomorphic to a subtournament of $G$,
and otherwise $G$ is {\em $H$-free}.  For a set $\mathcal{H}$ of tournaments,
$G$ is {\em $\mathcal{H}$-free} if $G$ is $H$-free for every 
$H \in \mathcal{H}$.
A set $X \subseteq V(G)$ is {\em transitive} if $G|X$ has no directed
cycles. The {\em chromatic number} of $G$ is the smallest integer $k$
for which $V(G)$ can be partitioned into $k$ transitive subsets.
Given tournaments $H_1$ and $H_2$
with disjoint vertex sets, we write $H_1 \Rightarrow H_2$ to mean the 
tournament $H$ with $V(H)=V(H_1) \cup V(H_2)$, and such that
$H|V(H_i)=H_i$ for $i=1,2$, and $V(H_1)$ is complete to $V(H_2)$.

A tournament $H$ is a {\em hero} if there exists $c$
(depending on $H$) such that every $H$-free tournament has chromatic
number at most $c$.  One of the results of \cite{heroes} is a 
complete characterization of all heroes. An important and the most difficult  
step toward that is the following:

\begin{theorem}
\label{chain}
If $H_1$ and $H_2$ are heroes, then so is $H_1 \Rightarrow H_2$.
\end{theorem}

It turns out that translating the proof of \ref{twographs} into the
language of tournaments gives a proof of \ref{chain} that is much simpler than
the one in \cite{heroes}, and we include  it here.

Let us start with some definitions. Let $H^*$ be  a tournament, and 
let $S$, $T$ be  subtournaments of $H^*$. As with undirected graphs,
we say that  $S$ is an {\em $H^*$-extension} of $T$ if $T=S \setminus u$ for 
some  $u \in V(S)$.
For an integer $k>0$ and  a set $\mathcal{F}$ of tournaments, we 
say that a tournament $G$ admits an $(\mathcal{F},k)$-partition if
$V(G) = X_1 \cup \ldots \cup X_k$, where for all $i \in \{1, \ldots, k\}$,   
either $|X_i|=1$, or  $G|X_i$ is $\mathcal{f}$-free. Please note that the
condition $|X_i|=1$ is significant only when all members of $\mathcal{F}$ 
have at most one vertex.

First we prove the tournament analogue of \ref{twographs}.

\begin{theorem}
\label{twotourn}
Let $H_1,H_2$ be non-null tournaments, and let $H$ be $H_1 \Rightarrow H_2$.
Let $m=\max(|V(H_1)|, |V(H_2)|)$, and $\mathcal{F} = \{H_1,H_2\}$. Then every $H$-free tournament $G$ admits an 
$(\mathcal{F}, 2(m+1)^m)$-partition.
\end{theorem}

\Proof 
By reversing all edges of $T$, if necessary, we may assume that $|V(H_1)|=m$.
We may assume that $G$ is not $H_1$-free, for otherwise the theorem holds.
Choose $h^*\in V(H_2)$, and let $H^*$ be the subtournament of $H$ with vertex set $V(H_1)\cup \{h^*\}$.
For a subtournament  $T$ of $H^*$, a  {\em $T$-piece} is an isomorphism
$g$ from $T$ to some subtournament of $G$ that we denote by $g(T)$.
For a $T$-piece $g$, and an
$H^*$-extension $S$ of $T$, let $V(S)\setminus V(T) = \{s\}$; we say that $v \in V(G) \setminus V(g(T))$ {\em $g$-corresponds to} $S$ if
the map sending $x$ to $g(x)$ for each $x\in V(T)$, and sending $s$ to $v$, is an $S$-piece.
We denote by $Y_S(g)$ the set of all vertices in $G$ that $g$-correspond to $S$.
Let $Y(g)$ be
the union of the sets $Y_S(g)$. Thus $Y(g)$ is the set of all
vertices in $V(G)$ that $g$-correspond to an $H^*$-extension of
$T$.
\\
\\
(1) {\em If $g$ is an $H_1$-piece, then
$Y(g)$ is complete from $V(g(H_1))$, and the graph $G|Y(g)$ is $H_2$-free.}
\\
\\
Since $V(H_1)$ is complete to $V(H_2)$ in $H$, and in particular complete to $h^*$, it follows that
$Y(g)$ is complete from $V(g(H_1))$. But then, since
$G$ is $H$-free, it follows
that $G|Y(g)$ is $H_2$-free. This proves~(1).
\\
\\
Let
$$\phi(t)= 
\begin{cases} 2(m+1)^{m-t} &\mbox{for } 0 \leq t \leq m-1\\
              1 &\mbox{for } t = m.
\end{cases}$$
\\
\\
(2) {\em Let $T$ be a subtournament of $H_1$, and let $g$
be a $T$-piece. Write $t=|V(T)|$. Then for some $H^*$-extension $S$ of $T$,
the tournament  $G|Y_S(g)$ admits an $(\mathcal{F}, \phi(t))$-partition.}
\\
\\
The proof is by induction on $m-t$.
If $t=m$, then (as $|H_1|=m$) it follows from~(1) that $G|Y(g)$ admits an $(\mathcal{F}, 1)$-partition,
so we may assume that $t \leq m-1$.

Choose an $H^*$-extension $S$ of $T$, with $S$ a subtournament of $H_1$.
If $G|Y_S(g)$ admits an $(\mathcal{F}, \phi(t))$-partition then we are done;
so we may assume that for every partition
of $Y_S(g)$ into at most $\phi(t)$ classes, some class contains both of
$H_1$, $H_2$.

Let $V(S)\setminus V(T) = \{s\}$.
Let $v \in Y_S(g)$, and let $h$ be the $S$-piece mapping $s$ to $v$ and mapping $x$ to $g(x)$ for each $x\in V(T)$.
Inductively, there exists an $H^*$-extension $Q$ of $S$
such that $G|Y_Q(h)$ admits an
$(\mathcal{F}, \phi(t+1))$-partition.
Let $V(Q)\setminus V(S) = \{r\}$, and let $R = Q\setminus s$.
We say that $v$ is {\em of type $R$}.
From the definition of $Y_Q(h)$, it follows that
$Y_Q(h) \subseteq Y_R(g)$,
and either
\begin{itemize}
\item $r$ is adjacent to $s$ in $H$, and $Y_Q(h)$ is the set of vertices in $Y_R(g)$ that
are different from and adjacent from $v$,
or
\item $s$ is adjacent to $r$ in $H$, and $Y_Q(h)$ is the set of
vertices in $Y_R(g)$ that are different from and adjacent to $v$.
\end{itemize}

For each $H^*$-extension $R$ of $T$ different from $S$, let $Z_{R}$ be the set of vertices in $Y_S(g)$ that are of type $R$.
Thus the sets $Z_{R}$ have union $Y_S(g)$.
Since $|V(H^*)|-|V(T)| = m+1-t$ and $R$ must be different from $S$, it follows that
are at most $m-t$ different types of vertices in $Y_S(g)$.
Since
$$m-t \leq 2(m+1)^{m-t}=\phi(t),$$
there is an $H^*$-extension $R$ of $T$ different from $S$, such that
$G|Z_{R}$  contains both of
$H_1$, $H_2$. Write $A=Y_R(g)$.
Let $V(R)\setminus V(T) = \{r\}$.

Assume first that $r$ is adjacent from $s$,
and choose $B \subseteq Z_{R}$ such that
$G|B$ is isomorphic to $H_2$. For each $b \in B$, let $A_b$ be the set
of vertices of $A \setminus B$ that are adjacent from $b$.
Let $A_0$ be the set of vertices of $A \setminus B$ that are complete
to $B$. Then $A \subseteq B \cup A_0 \cup \bigcup_{b \in B}A_b$.

Let $b \in B$; since $b \in Z_{R}$, since
$A_b \subseteq Y_R(g)$, and since $A_b$ is complete from $b$,
it follows from the definition of $Z_{R}$ that $G|A_b$ admits an
$(\mathcal{F}, \phi(t+1))$-partition.
Since $G$ is $H$-free, it follows that $G|A_0$ is $H_1$-free.
Since $|B| \leq m$, this  implies that $G|A$ admits an
$(\mathcal{F},K)$-partition, where
$$K= m \phi(t+1)+ m +1 \leq 2(m+1)^{m-t}=\phi(t),$$
and so (2) holds.

Now we assume that $r$ is adjacent to $s$.
Choose $B \subseteq Z_{R}$ such that
$G|B$ is isomorphic to $H_1$. For $b \in B$, let $A_b$ be the set
of vertices of $A \setminus B$ that are adjacent to $b$.
Let $A_0$ be the set of vertices of $A \setminus B$ that are complete from 
$B$. Then $A \subseteq B \cup A_0 \cup \bigcup_{b \in B}A_b$.

Let $b \in B$;  since $b \in Z_{R}$, since
$A_b \subseteq Y_R(g)$, and since $A_b$ is complete to $b$,
it follows from the definition of $Z_{R}$ that $G|A_b$ admits an
$(\mathcal{F}, \phi(t+1))$-partition.
Since $G$ is $H$-free, it follows that $G|A_0$ is $H_2$-free.
Since $|B| \leq m$, this  implies that $G|A$ admits an
$(\mathcal{F},K)$-partition, where
$$K = m \phi(t+1) + m +1 \leq 2(m+1)^{m-t}=\phi(t).$$
This proves~(2).
\\
\\
Now let $T$ be the null tournament, and $g$ the isomorphism from $T$ into $G$. Then $T$ is an
subtournament of $H_1$, and $g$ is a $T$-piece. Also,
$Y_S(g) = V(G)$ for every $H^*$-extension $S$ of $T$.
But then, by~(2), $G$ admits an
$(\mathcal{F}, 2(m+1)^m)$-partition. This proves~\ref{twotourn}.~\bbox

Now \ref{chain} follows easily:

{\bf Proof of \ref{chain}.} Since $H_1$ and $H_2$ are heroes, there
exists an integer $c>0$ such that every $H_i$-free tournament has
chromatic number at most $c$ for $i=1,2$. By \ref{twotourn}, 
every $H$-free tournament $G$ has an $(\{H_1,H_2\}, 2(m+1)^m)$-partition,
where  $m=\max(|V(H_1)|, |V(H_2)|)$; and therefore $V(G)$ can be
partitioned into $2(m+1)^mc$ transitive subsets. Thus every
$H$-free tournament has chromatic number at most  $2(m+1)^mc$,
and consequently $H$ is a hero. This proves \ref{chain}.~\bbox

\section{Cographs}

In this section we prove \ref{cographsplit}, which is
the cograph analogue of \ref{Ksplit}.
Let $\mathcal{F}$ be a set of graphs,
where $k \geq 1$ is an integer, and let $P>0$ be an integer. 
We say that a graph $C$ is  {\em $(\mathcal{F},P)$-universal} if for
every partition $X_1, \ldots, X_P$ of $V(C)$ (by a {\em partition} we
mean that the sets  $X_1, \ldots, X_P$ are pairwise disjoint,
and have union $V(C)$), there exists $i \in \{1, \ldots, P\}$ such that
$C|X_i$ contains every member of $\mathcal{F}$ (in other words, $C$ does
not admit an $(\mathcal{F},P)$-partition.

We start with a lemma that establishes the existence of universal cographs.
\begin{theorem}
\label{universal}
Let $P,k$ be positive integers.
Let $\mathcal{F}$ be a set of connected cographs, all of height at most $k$. 
Then there exists a connected cograph of height $k$ that is 
$(\mathcal{F},P)$-universal.
\end{theorem}

\Proof The proof is by induction on $k$. Suppose first that  
$k=1$. Then  the members  of $\mathcal{F}$ are complete graphs; let 
$m=\max_{F \in \mathcal{F}}|V(F)|$. Now the complete  graph on
$mP$ vertices is $(\mathcal{F},P)$-universal.

Next we consider a general $k>1$. For every $F \in \mathcal{F}$,
the members of the set $ac(F)$ are anticonnected cographs of
height at most $k-1$. Let $A=\bigcup_{F \in \mathcal{F}}ac(F)$. Inductively, 
passing to the complement, there exists an $(A,P)$-universal anticonnected
cograph $C$ of height $k-1$.  Denote by $s$ the maximum number
of anticomponents of a member of $\mathcal{F}$, and write $K=(s-1)P+2$.
Then $K \geq 2$.
Let $U$ be the cograph obtained from $K$ vertex-disjoint copies 
$C_1,\ldots, C_K$ of $C$ by making $V(C_i)$ complete to $V(C_j)$ for 
all $1 \leq i<j \leq K$. Then $U$ is a connected cograph of height $k$.

We claim that $U$ is $(\mathcal{F},P)$-universal. Let $X_1, \ldots, X_P$
be a partition of $V(U)$. We need to prove that $U|X_j$ 
contains every member of $\mathcal{F}$ for some $j \in \{1, \ldots, P\}$.

For $i \in \{1, \ldots, K\}$ and 
$j \in \{1, \ldots, P\}$ write $C_i^j=C_i\cap X_j$.
Since $C$ is $(A,P)$-universal, it follows that for every 
$i \in \{1, \ldots, K\}$ there exists $j \in \{1, \ldots, P\}$
such that $C_i^j$ contains every member of $A$. 
For every $j \in \{1, \ldots, P\}$, let 

$$I_j=\{i \in \{1, \ldots, K\} \text{ such that } C_i^j \text{ contains every 
member of } A \}.$$

Since $K>P(s-1)$, there exists $j \in \{1, \ldots, P\}$ such that 
$|I_j| \geq s$. Let $D$ be the graph obtained from the graphs
$\{C_i^j\}_{i \in I_j}$   by making $V(C_i^j)$ complete to 
$V(C_h^j)$ for all distinct $i,h \in I_j$. Then  each anticomponent of $D$ 
contains every member of $A$, and $D$ has at least $s$ anticomponents. But 
since each member of $\mathcal{F}$ has at most $s$
anticomponents, it follows that $D$ contains every member of $\mathcal{F}$.
Since $D$ is an induced subgraph of $U|X_j$, it follows that 
$U|X_j$ contains every member of $\mathcal{F}$.
This proves the claim that $U$ is  $(\mathcal{F},P)$-universal, and
completes the proof of \ref{universal}.~\bbox

We are now ready to prove~\ref{cographsplit}. First we note that
\ref{disconnected} immediately implies the following:

\begin{theorem}
\label{cographs2}
Let $k \geq 0$ be an integer, and let  $H$ and $J$ be cographs, each of height
$k+1$, such that $H$ is anticonnected, and $J$ is connected. Then there exists
an integer $P$ such that every $\{H,J\}$-free  graph admits a 
$(c(H) \cup ac(J), P)$-partition.
\end{theorem}

\ref{cographsplit} says:

\begin{theorem}
\label{cographsplit2}
Let $k \geq 1$ be an integer, and let  $H$ and $J$ be cographs, each of height
$k+1$, such that $H$ is anticonnected, and $J$ is connected. Then there
exist cographs $\tilde{H}$ and $\tilde{J}$, each of height $k$, such that
 $\tilde{H}$ is connected, and $\tilde{J}$ is anticonnected, and 
every $\{H,J\}$-free  graph $G$ is $\{\tilde{H},\tilde{J}\}$-split.
\end{theorem}

\Proof Write $A=c(H)$ and $B=ac(J)$. Then the members of $A$ are connected
cographs of height at most $k$, and the members of $B$ are anticonnected 
cographs of height at most $k$. By~\ref{cographs2}, $G$  admits an 
$(A \cup B, P)$-partition. We may assume that there exists 
$j  \in \{1, \ldots, P\}$ such that for $1 \leq i \leq j$
either $|X_i|=1$, or 
the subgraph of $G$ induced by $X_i$ is 
$\{F\}$-free  for some $F \in A$, and for  $j < i \leq P$
either $|X_i|=1$, or 
the subgraph of $G$ induced by $X_i$ is 
$\{F\}$-free  for some $F \in B$.
Let $X=\bigcup_{1 \leq i \leq j}X_i$,
and $Y=\bigcup_{j<i \leq P}X_i$. By~\ref{universal}, there exists 
a connected cograph $\tilde{H}$, of height $k$, such that   $\tilde{H}$ is 
$(A,P)$-universal. By~\ref{universal} (complemented) there 
exists an anticonnected cograph $\tilde{J}$, of height $k$, such that  
$\tilde{J}$ is  $(B,P)$-universal. From the definition of $X$ and $Y$, it 
follows that $G|X$ is $\tilde{H}$-free, and $G|Y$ is $\tilde{J}$-free, and so
$G$ is $\{\tilde{H},\tilde{J}\}$-split, as required. This 
proves~\ref{cographsplit2}.~\bbox

\section{A construction}

Let $G$ be a graph. A {\em block} of $G$ is a maximal subgraph 
of $G$ that is either $2$-connected 
or isomorphic to $K_2$ (so in this paper isolated vertices do not belong to a block). 
Please note that two distinct blocks cannot share more than one vertex.
In this section we prove \ref{construction}, which we restate:

\begin{theorem}
\label{construction2}
Let $L$, $M$ be graphs each with at least one edge.  Then for every choice of 
non-negative integers $r,k$ there is a graph $G$ such that 
\begin{itemize}
\item for every $S \subseteq V(G)$ with $|S| \leq r$, the graph
$G|S$ is $\{L,M\}$-split, and
\item $G$ has no $(\{L,M\}, k)$-partition.
\end{itemize}
\end{theorem}

\Proof 
Let $L_1,\dots,L_l$ be the blocks of $L$ and let 
$M_1,\dots,M_m$ be the blocks of $M$.  We may assume that $L_1,\dots,L_l,M_1,\dots,M_m$ each have at most $|V(L_1)|$ vertices, 
and $l,m\geq1$ (since $L,M$ each have at least one edge by hypothesis). Note that $L,M$ may not be connected, and there may be isolated vertices that are not contained in any block.
Fix $r,k$; we may assume that $r\ge\max\{ 3|L|,3|M|\}$. Choose a small constant $\epsilon>0$ ($\epsilon=1/(r+2)$ will do).

Let $V$ be a set of size $n\ge1$.
By a {\em hypergraph} with vertex set $V$ we mean in this paper a set of subsets of $V$ (all different), and we call these subsets {\em hyperedges}.
We generate $l+m$ independent random hypergraphs with vertex set $V$ as follows.  For $i=1,\dots,l$, we let $H_i^L$ be a random $|V(L_i)|$-uniform 
hypergraph, where each possible hyperedge is present independently with probability $p_i=n^{-(|V(L_i)|-1)+\epsilon}$; 
for $j=1,\dots,m$, we let $H_j^M$ be a random $|V(M_j)|$-uniform hypergraph, where each possible hyperedge is present independently 
with probability $q_j=n^{-(|V(M_j)|-1)+\epsilon}$.  We then set $H$ to be the union of these $l+m$ hypergraphs, labeling each hyperedge of $H$ with 
the name of the hypergraph ($L_i$ or $M_j$) it came from; we refer to the $L_i$ and $M_j$ as {\em pieces of $H$}.  Note that at this point 
a hyperedge might have more than one label. 

For $t\ge3$, a {\em cycle of length $t$} in $H$ (or {\em $t$-cycle}) is a sequence of distinct vertices $v_1,\dots,v_t$ such that for $1\le i\le t$, some hyperedge $A$
satisfies 
$A\cap\{v_1,\ldots,v_t\}=\{v_i,v_{i+1}\}$, 
where subscripts are taken modulo $t$. A {\em cycle of length $2$} in $H$ (or {\em $2$-cycle}) is a pair of distinct vertices $v_1,v_2$ such that there
are at least two distinct hyperedges containing $v_1,v_2$.
If $R\subseteq V$, we denote by $H\setminus R$ the hypergraph of all hyperedges of $H$ that are disjoint from $R$.
\\
\\
(1) {\em If $n$ is sufficiently large, then with high probability, there is a set $R$ of size $o(n/\log n)$ such that $H\setminus R$ has
no cycles of length at most $r$ and no hyperedges with multiple labels.}
\\
\\
Let us deal first with multiple labels.  If a hyperedge of size $k$ has at least two labels, then it has been chosen in (at least) two distinct pieces of $H$.  
This has probability at most $\binom{l+m}{2}(n^{-(k-1)+\epsilon})^2=O(n^{-2(k-1)+2\epsilon})$.  
The expected number of such hyperedges is $O(n^{-2(k-1)+2\epsilon}\binom nk)=O(n^{-k+2+2\epsilon})=O(n^{2\epsilon})$, so by Markov's Inequality 
there are with high probability at most $o(n/\log n)$ such hyperedges.  

For vertices $v,w$ in $H$, let $X_{vw}$ be the number of hyperedges of $H$ that contain both $v$ and $w$.  Then
\begin{align*}
\E X_{vw}
&=\sum_i \binom{n-2}{|V(L_i)|-2}p_i+\sum_j \binom{n-2}{|V(M_j)|-2}q_j\\
&\le \sum_i n^{|V(L_i)|-2} n^{-(|V(L_i)|-1)+\epsilon} + \sum_j n^{|V(M_j)|-2} n^{-(|V(M_j)|-1)+\epsilon}\\
&=O(n^{-1+\epsilon}).
\end{align*}
The probability that there is a pair of hyperedges both containing the pair $\{v,w\}$ is at most 
$\E\binom{X_{vw}}{2}$,  
the expected
number of pairs of hyperedges containing $\{v,w\}$.  Since $X_{vw}$ is a sum of independent indicator variables, we have
$\E X_{vw}(X_{vw}-1) \le (\E X_{vw})^2$ and so
$\E\binom{X_{vw}}{2}=O(n^{-2+2\epsilon})$.  Summing over all $v,w$, we see that the expected number of 
pairs $\{v,w\}$ that lie in two or more hyperedges is $O(n^{2\epsilon})$ and so by Markov's Inequality is with high probability $O(n^{3\epsilon})$.
Consequently with high probability the number of $2$-cycles is $O(n^{3\epsilon})$.

For fixed $t\ge3$, we now bound the number of $t$-cycles.  Let $v_1,\dots,v_t$  be a sequence of distinct vertices, and for $i=1,\ldots t$, 
let $Y_i$ be the number of hyperedges that meet $\{v_1,\ldots,v_t\}$ in exactly $\{v_i,v_{i+1}\}$ (subscripts taken modulo $t$).  
Note that the random variables $Y_1,\dots,Y_t$ are independent, as they depend on disjoint sets of hyperedges; also, for each $i$, 
we have $Y_i\le X_{v_iv_{i+1}}$.  Then the probability that the sequence $v_1,\dots,v_t$ forms a $t$-cycle is
$$
\P[Y_1\cdots Y_t>0]\le\E[Y_1\cdots Y_t]=\prod_{i=1}^t \E Y_i\le\prod_{i=1}^t \E X_{v_iv_{i+1}} 
=O( n^{-t+t\epsilon}).
$$
Summing over all $O(n^t)$ choices of  $v_1,\dots,v_t$, we see that the expected number of $t$-cycles is
$O( n^{t\epsilon})$.  
So by Markov's Inequality, with high probability the number of $t$-cycles is $O(n^{(t+1)\epsilon})=O(n^{(r+1)\epsilon})$.
Consequently, with high probability the number of cycles of length at most $r$ is at most $O(n^{(r+1)\epsilon})=o(n/\log n)$.

Finally, let $R$ consist of one vertex from each hyperedge with multiple labels, and one vertex from each cycle of length at most $r$.
By the argument above, this gives with high probability a total of at most $o(n/\log n)$ vertices. This proves~(1).
\\
\\
(2) {\em  If $n$ is sufficiently large, then with high probability, every set of at least $n/\log n$ vertices contains 
hyperedges from every $H_i^L$ and $H_j^M$.}
\\
\\
Let $S$ be a set of at least $n/\log n$ vertices.  Then, for any $i$, the probability that $H$ contains no hyperedge from $H_i^L$ is at most
$$(1-p_i)^{\binom{|S|}{|V(L_i)|}}
\le (1-p_i)^{n^{|V(L_i)|-\epsilon/2}}
\le \exp(-p_in^{|V(L_i)|-\epsilon/2})
= \exp(-n^{1+\epsilon/2}),
$$
provided $n$ is sufficiently large.  The same bound holds for hyperedges from $H_j^M$.  There are fewer than $2^n$ choices for $S$, 
and $O(1)$ choices of $i$ or $j$, so with high probability every set of at least $n/\log n$ vertices contains hyperedges from every $H_i^L$ and $H_j^M$. This proves~(2).
\\
\\
From (1) and (2), if $n$ is sufficiently large then there exists a hypergraph $H$ and a subset $R\subseteq V$ such that
\\
\\
(3) {\em $R$ has size $o(n/\log n)$, and $H\setminus R$ has
no cycles of length at most $r$, and has no hyperedges with multiple labels, and every set of at least $n/\log n$ vertices of $V$ contains 
hyperedges of $H$ from every $H_i^L$ and $H_j^M$.}
\\
\\
Let $H'=H\setminus R$.  
We construct a graph $G$ by replacing each surviving hyperedge of form $H_i^L$ by a copy of $L_i$, and each surviving hyperedge 
of form $H_i^M$ by a copy of $M_i$ (in each case, choosing an arbitrary ordering of the vertices).  
Note that this is well-defined, as $H'$ has no hyperedges with multiple labels, and no pair of vertices belongs to two hyperedges 
(so the subgraphs we are inserting intersect pairwise in at most one vertex).

We claim that, provided $n$ is sufficiently large, $G$ satisfies the theorem.
\\
\\
(4) {\em  For every subset $S$ of $V(G)$ with size at most $r$,
the graph $G|S$ is $(L,M)$-split.} 
\\
\\
Fix an $S$.  We show that $S$ can be partitioned into two sets, $X$ and $Y$,
so that $G|X$ is $L$-free, and $G|Y$  is $M$-free. Let $H_S$ be the hypergraph
containing all sets $A\cap S$, where $A\in H'$ and $|A\cap S|\ge 2$.
We label each hyperedge of $H_S$ with the label of the hyperedge that generated it 
(note that this is well defined: all hyperedges of $H_S$ 
have size at least two, and so are contained in only one hyperedge of $H'$).

By construction, the hypergraph $H_S$ has no cycles, since any cycle in $H_S$ 
is a cycle in $H'$, and $H'$ has no cycles of length at most $r$.  
It is straightforward to find a partition $(X,Y)$ of $S$ such that 
every hyperedge of $H_S$ has exactly one vertex in $Y$.

Suppose that $G|X$ is not $L$-free. Then, in particular, there exists a 
subset $B$ of $X$ such that $G|B$ is  isomorphic to $L_1$.
Since $L_1$ is 2-connected, and $H_S$ has no cycles, it follows that  $B$ is 
contained in some hyperedge $E$ of $H_S$. But all hyperedges of $H_S$ have size at
most $|V(L_1)|$, and therefore $E\cap Y = \emptyset$, a contradiction. 
This proves that $G|X$ is  $L$-free.

Next suppose that $G|Y$ is not $M$-free. Then, in particular, there exists a 
subset $B$ of $Y$ such that $G|B$ is  isomorphic to $M_1$.
Since $M_1$ is 2-connected, and $H_S$ has no cycles, it follows that
$B$ is contained  in some hyperedge $E$ of $H_S$.  But $|E\cap Y| = 1$, a contradiction since $|V(M_1)| \geq 2$. Thus
$G|Y$ is $M$-free. This proves~(4).
\\
\\
(5) {\em $G$ has no $(\{L,M\},k)$-partition.}
\\
\\
It is enough to show that for every subset $S$ of $V(G)$ with $|S|\ge n/2k$,
the graph $G|S$ is not $L$-free and not $M$-free. 
Let ${\mathcal B}=\{L_1,\dots,L_l,M_1,\dots,M_m\}$.
Note that, by adding fewer than $2|L|$ additional blocks each isomorphic to $L_1$, we can generate a connected graph $L'$ that has $L$ as an induced subgraph, 
and such that all its blocks belong to ${\mathcal B}$; and similarly for $M$.  
Thus, since $r\ge \max\{3|L|, 3|M|\}$,  it will be enough to prove the following.
\\
\\
(6) {\em For every integer $t\ge1$ with $t<r$, there is an integer $K(t) \geq 0$ such that if $n$ is sufficiently large then the following holds.
Let $F$ be a connected graph with exactly $t$ blocks, such that all its blocks are isomorphic to members of ${\mathcal B}$; then 
for every set $W \subseteq V(G)$ of at least $K(t)n/\log n$ vertices, there is an induced subgraph $F'$ of $G|W$, isomorphic to $F$,
such that the vertex set of every block of $F'$ is a hyperedge of $H'$.}
\\
\\
We argue by induction on $t$.  For $t=1$ this follows from (3), taking $K(1) = 1$.
So suppose $t\ge2$ and we have shown the existence of $K(t-1)$.  Write 
$K=K(t-1)$.
Let $F$ be a graph with $t$ blocks, all from 
$\mathcal B$.   Since $t \geq 2$ and $F$ is connected, and the bipartite graph of blocks versus cutpoints of $F$ is a tree,
it follows that there is a block $B$ of $F$, and a vertex $v_0$ of $B$, such that no other block of $F$ contains any vertex of $B$
different from $v_0$.
Write $F'=F\setminus (V(B)\setminus \{v_0\})$. Then $F'$ has $t-1$ blocks.
If $P$ is an induced subgraph of $G$ and $v\in V(P)$, and there is an isomorphism between $P$ and $F'$
mapping $v$ to $v_0$, we call $v$ an {\em anchor} of $P$.

Pick a large constant $M$, let $W$ be any set of at least $Mn/\log n$ vertices in $G$, 
and let ${\mathcal F}'$ be a maximal collection of pairwise vertex-disjoint copies of $F'$ with vertices from $W$.  
Then, by our inductive hypothesis, the union
of the vertex sets of the members of ${\mathcal F}'$ contains all but at most 
$Kn/\log n$ vertices from $W$, and so 
$$|{\mathcal F}'|\ge \frac{(M-K)n}{|V(F')|\log n}>\frac{n}{\log n},$$
provided $M$ is sufficiently large.
Let $T$ be a set consisting of an anchor of each member of  $\mathcal F'$.
Then $|T|\ge n/\log n$, and since $K(1) = 1$, $T$ contains an
hyperedge $E$ of $H$ such that $G|E$ is isomorphic to $B$. Let $z\in E$ such that some such isomorphism takes $z$ to $v_0$.
Let $P \in \mathcal{F}'$ be such that $E \cap V(P) =\{z\}$.
If some vertex in $E\setminus \{z\}$ is adjacent in $G$ to some vertex in $V(F)\setminus \{z\}$, then since every edge of $G$ is contained
in a hyperedge of $H'$, and $E$ is a hyperedge of $H'$, and $P$ has at most $r-2$ blocks, each with vertex set some hyperedge of $H'$, it follows that
$H'$ has a cycle of length at most $r$, a contradiction.
Thus there is no such edge, and so $G|(V(P) \cup E)$ is isomorphic to $F$, and (6) holds taking $K(t)=M$.
This proves~(6) and completes the proof of~\ref{construction2}.~\bbox

\end{document}